\documentclass[review]{elsarticle}
\usepackage{algorithm, algorithmic}
\usepackage{lineno,hyperref}
\usepackage{mathdots}

\modulolinenumbers[5]

\journal{Journal of \LaTeX\ Templates}

\begin{document}

\begin{frontmatter}

\title{A Breakdown Free Numerical Algorithm for Inverting General Tridiagonal Matrices}
\tnotetext[mytitlenote]{Fully documented templates are available in the elsarticle package on \href{http://www.ctan.org/tex-archive/macros/latex/contrib/elsarticle}{CTAN}.}

%% Group authors per affiliation:
\author{Moawwad El-Mikkawy\corref{mycorrespondingauthor}}
\author{Abdelrahman Karawia}
\address{Mathematics Department, Faculty of Science, Mansoura University, Mansoura, 35516, Egypt}
\cortext[mycorrespondingauthor]{Corresponding author}

\begin{abstract}
In the current paper the authors linked two methods in order to evaluate general n-th order tridiagonal determinants. A breakdown free numerical algorithm is developed for computing the inverse of any $n\times n$ general nonsingular tridiagonal matrix without imposing any constrains. The algorithm is suited for implementation using any computer language such as FORTRAN, PYTHON, MATLAB, MAPLE, C, C++, MACSYMA, ALGOL, PASCAL and JAVA. Some illustrative examples are presented.
\end{abstract}

\begin{keyword}
Tridiagonal matrix\sep Determinants\sep  Matrix inverse\sep Numerical algorithm\sep Hadamard product 
\MSC[2010] 15A09\sep 15A15\sep 39A06\sep 65F05\sep 65Y20\sep  68W30
\end{keyword}

\end{frontmatter}

\linenumbers

\section{Introduction}
Applications of matrices in real life include physics, geology, statistics, graphics, economic, business, animation, encryption, games and construction.\\
\textbf{Definition 1.} A band matrix $A_n=(t_{ij})_{i,j=1}^n$ is called tridiagonal matrix if it takes the form:\\
\begin{equation}
    A_n=[t_{ij}]_{i,j=1}^n=trid(b_i,d_i, a_i)= \left[
\begin{array}{cccccc}
    d_1 & a_1 & 0 & \cdots & \cdots & 0\\
   b_1 & d_2 & a_2 &\ddots &  & \vdots \\
   0 & b_2 & d_3 & \ddots & 0 & \vdots\\
   \vdots & \ddots &\ddots &\ddots &\ddots & 0\\
   \vdots & & 0 &\ddots &\ddots & a_{n-1}\\ 
   0& \cdots & \cdots & 0 & b_{n-1} & d_n\\
 \end{array}
\right]
\end{equation}
in which $t_{ij}=0$ whenever $|i-j| > 1$.\\
Band matrices of tridiagonal type frequently appear in scientific and engineering areas (see for instance \cite{USMANI1994a,Usmani1994,Hadj2008,Ran2009,RAN2009b}). The need to obtain the inverse of nonsingular tridiagonal matrices arises in such cases. This problem has been investigated by many authors. The interested reader may refer to \cite{Caratelli2021,Kumar1993,Huang1997,Li2010,ABDERRAMANMARRERO2013,Elmikkawy2004} and the references therein.\\
Perhaps the first complete analysis for inverting nonsingular tridiagonal matrices without imposing any constraints was considered in \cite{Huang1997}. Nevertheless the resulting numerical algorithm is not breakdown free. It breaks down if any principal submatrix is singular. By introducing the symbolic computation technology, we analyzed and established an algorithm which is breakdown free \cite{Elmikkawy2006}. The computational complexity of the algorithms given in \cite{Huang1997,Elmikkawy2006} is $O(n^2)$.\\
In the applied domain, the numerical approach is generally preferable \cite{ABDERRAMANMARRERO2013}. This is in fact the motivation of the present paper. We are mainly concerned with constructing efficient numerical algorithms for computing the determinant and the inverse of a general n-th order tridiagonal
matrix.\\
The present article is organized as follows. A cost-efficient numerical algorithm for evaluating general n-th order tridiagonal determinants is presented in Section 2. In Section 3, the main result is given. Some illustrative examples are also introduced.

\section{A Cost-efficient Numerical Algorithm for Evaluating General n-th Order Tridiagonal Determinants}
In this section, we focus on constructing an efficient numerical algorithm for the determinant of $A_n$ in (1). The algorithm is mainly based on the following facts:
\begin{enumerate}[(1)]
    \item The reversal matrix $J_n$ defined by:\\
    \begin{equation}
        J_n=\left[
\begin{array}{cccccc}
    0 & 0 & \cdots & 0 & 0 & 1\\
   0 & 0 &\cdots & 0 & 1 & 0 \\
   0 & 0 & \cdots & 1 & 0 & 0\\
   \vdots & \vdots &\iddots &\vdots &\vdots & \vdots\\
   0 & 1 & \cdots &\cdots & \cdots\\ 
   1& 0 & \cdots & 0 & 0 & 0\\
 \end{array}
\right]
    \end{equation}
    satisfies:\\
    (i) $det(J_n)=\pm 1$, and\\
    (ii) $J_n=J_n^t=J_n^{-1}$.
    \item Let $M=[m_{ij}]$ be an $n\times n$ matrix. The rotate of $M$, denoted $M^R$ is defined as: $M^R=J_nMJ_n$. If $M=M^R$, then $M$ is called centrosymmetric matrix. It is true that $det(M)=det(M^R)$.
\end{enumerate}
To study tridiagonal matrices, it is advantageous to introduce two n-dimensional vectors, $\mathbf{c}=[c_1, c_2, \cdots, c_n]^t$, 
and $\mathbf{e}=[e_1, e_2, \cdots, e_n]^t$, and two (n+1)-dimensional vectors, $\mathbf{f}=[f_0,f_1, f_2, \cdots, f_n]^t$, and $\mathbf{g}=[g_1, g_2, \cdots, g_n, g_{n+1}]^t$. These vectors are given by:\\
\begin{equation}
    c_1=d_1, c_i=d_i-\frac{a_{i-1}b_{i-1}}{c_{i-1}}, i=2, 3, \cdots, n,
\end{equation}
\begin{equation}
    e_n=d_n, e_i=d_i-\frac{a_ib_i}{e_{i+1}}, i=n-1, n-2, \cdots, 1,
\end{equation}
\begin{equation}
    f_0=1, f_i=\left|
\begin{array}{cccccc}
    d_1 & a_1 & 0 & \cdots & \cdots & 0\\
   b_1 & d_2 & a_2 &\ddots &  & \vdots \\
   0 & b_2 & d_3 & a_3 & 0 & \vdots\\
   \vdots & \ddots &\ddots &\ddots &\ddots & 0\\
   \vdots & & 0 &\ddots &\ddots & a_{i-1}\\ 
   0& \cdots & \cdots & 0 & b_{i-1} & d_i\\
 \end{array}
\right|, i=1, 2, 3, \cdots, n,
\end{equation}
and
\begin{equation}
  g_{n+1}=1, g_i=\left|
\begin{array}{cccccc}
    d_i & a_i & 0 & \cdots & \cdots & 0\\
   b_i & d_{i+1} & a_{i+1} &\ddots &  & \vdots \\
   0 & b_{i+1} & d_{i+2} & a_{i+2} & 0 & \vdots\\
   \vdots & \ddots &\ddots &\ddots &\ddots & 0\\
   \vdots & & 0 &\ddots &\ddots & a_{n-1}\\ 
   0& \cdots & \cdots & 0 & b_{n-1} & d_n\\
 \end{array}
\right|, i=n, n-1, n-2, \cdots, 1.  
\end{equation}
The four vectors are related by \cite{ELMIKKAWY2003}:
\begin{equation}
    f_i=\prod_{r=1}^ic_r=c_if_{i-1}, i=1, 2, \cdots, n,
\end{equation}
and
\begin{equation}
    g_i=\prod_{r=i}^ne_r=e_ig_{i+1}, i=n, n-1, \cdots, 1
\end{equation}
Consequently,
\begin{equation}
   f_n=det(A_n)=g_1 
\end{equation}
The principal minors $f_i$'s and $g_i$'s in (5) and (6) satisfy the second order linear difference equations:
\begin{equation}
    f_i=d_if_{i-1}-a_{i-1}b_{i-1}f_{i-2},\quad i=2, 3, \cdots, n, 
\end{equation}
\hspace{2.6cm}with $f_0=1, f_1=d_1$.\\
and
\begin{equation}
    g_i=d_ig_{i+1}-a_ib_ig_{i+2},\quad i=n-1, n-2, \cdots, 1.
\end{equation}
\hspace{2.6cm}with $g_{n+1}=1, g_n=d_n$\\
respectively.
Note that for each $i=1, 2, \cdots, n$, $f_i$ and $g_i$ are determinants of order $i$ and $n+1-i$, respectively. By using the syntax of Maple language, we have:
$$f_i=det(submatrix(A_n,1..i,1..i)) \quad \textit{and}\quad g_i=det(submatrix(A_n,i..n,i..n))$$
It is worth mentioned that any of the following conditions is sufficient (but not necessary) to ensure that the tridiagonal matrix $A_n$ in (1) is nonsingular.
\begin{itemize}
    \item (i) All $c_i\ne 0, i=1, 2, \cdots, n$.
    \item (ii) All $e_i\ne 0, i=1, 2, \cdots, n$.
\end{itemize}
For example, consider the two matrices $A$ and $B$ given by:
$$A=\left[
\begin{array}{ccc}
    1 & 1 & 0\\
   1 & 1 & 2\\
   0 & 2 & 3\\
 \end{array}
\right], B=\left[
\begin{array}{ccc}
    3 & 2 & 0\\
   -1 & 1 & 1\\
   0 & 1 & 1\\
 \end{array}
\right]$$\\
For these matrices, we have\\ $det(A)=-4\ne 0$ and $det(B)=2\ne 0$ although $c_2=e_2=0$.\\
Also, each of the following conditions is necessary and sufficient to ensure that a symmetric tridiagonal matrix of the form (1) is positive definite\cite{Elmikkawy2006}.
\begin{itemize}
    \item (i) All $c_i> 0$ for all $i=1, 2, \cdots, n$.
    \item (ii) All $e_i> 0$ for all $ i=1, 2, \cdots, n$.
\end{itemize}
The nonsingularity of a matrix should be checked before considering the problem of computing its inverse. So, we will present two numerical algorithms for evaluating the determinant of $A_n$. \\
For both \textbf{Algorithm 1} and \textbf{Algorithm 2}, the authors successfully linked two methods together in order to evaluate n-th order tridiagonal determinants. Throughout this paper $\sum_{k+1}^k (.)$ is assumed equal to $0$ and $\prod_{k+1}^k (.)$ is assumed equal to $1$. Also, $p_i=-a_i$ and $q_i=-b_i, i=1, 2, \cdots, n-1$. 
\\
\begin{algorithm}[!h]
 \caption{(To evaluate $f_1, f_2, \cdots, f_n=det(A_n)$)}
  \textbf{Input:} $n$ and the components of the vectors $\mathbf{a}, \mathbf{b}$, and $\mathbf{d}$.\\
  \textbf{Output:} $det(A_n)$.\\
  \textbf{Step 1:} Set $c_1=d_1$, $f_1=d_1$, and $m=1$.\\
  \textbf{Step 2:} While $m\le n-1$ and $c_m\ne 0$ do\\
  \hspace*{2cm}  $m=m+1$,\\
  \hspace*{2cm}  $c_m=d_m - a_{m-1}b_{m-1}/c_{m-1}$,\\
  \hspace*{2cm}  $f_m=c_m f_{m-1}$,\\
  \hspace*{1.5cm}End do.\\
  \textbf{Step 3:} For $k=m+1$ to $n$ do\\
  \hspace*{2cm}  $f_k=d_k f_{k-1}-a_{k-1} b_{k-1} f_{k-2}$\\
  \hspace*{1.5cm}End do.\\
  \textbf{Step 4:} Set $det(A_n)=f_n$.
  \label{alg1}
  \end{algorithm}
 
 \begin{algorithm}[!h]
  \caption{(To evaluate $g_n, g_{n-1}, \cdots, g_1=det(A_n)$)}
  \textbf{Input:} $n$ and the components of the vectors $\mathbf{a}, \mathbf{b}$, and $\mathbf{d}$.\\
  \textbf{Output:} $det(A_n)$.\\
  \textbf{Step 1:} Set $e_n=d_n$, $g_n=d_n$, and $m=n$.\\
  \textbf{Step 2:} While $m\ge 2$ and $e_m\ne 0$ do\\
  \hspace*{2cm}  $m=m-1$,\\
  \hspace*{2cm}  $e_m=d_m - a_mb_m/e_{m+1}$,\\
  \hspace*{2cm}  $g_m=e_m g_{m+1}$,\\
  \hspace*{1.5cm}End do.\\
  \textbf{Step 3:} For $k=m-1$ step $-1$ to $1$ do\\
  \hspace*{2cm} $g_k=d_k g_{k+1}-a_k b_k g_{k+2}$\\
  \hspace*{1.5cm}End do.\\
  \textbf{Step 4:} Set $det(A_n)=g_1$.
  \label{alg2}
  \end{algorithm}
 The \textbf{Algorithm 1} and \textbf{Algorithm 2} will be referred to as \textbf{DETGTRIF}  and \textbf{DETGTRIG} respectively.
\section{The Main Result}
This section is mainly devoted for developing a new efficient numerical algorithm for computing the inverse of general nonsingular tridiagonal matrices of the form (1).\\
\textbf{Lemma 1:} For $k=1, 2, \cdots, n$, we have
\begin{equation}
    f_kg_{k+1}-d_kf_{k-1}g_{k+1}+f_{k-1}g_k=f_n
\end{equation}
\textbf{Proof:} for $k=1$, the left-hand side of (12) gives\\
$L.H.S=f_1g_2-d_1f_0g_2+f_0g_1=(d_1)(g_2)-(d_1)(1)(g_2)+(1)g_1=g_1=f_n=R.H.S$\\
Thus formula (12) holds for $k=1$.\\
We next assume that formula (12) holds for $k=1, 2, \cdots, i-1$. That is
\begin{equation}
    f_{i-1}g_i-d_{i-1}f_{i-2}g_i+f_{i-2}g_{i-1}=f_n
\end{equation}
We want to prove that:
\begin{equation}
    f_ig_{i+1}-d_if_{i-1}g_{i+1}+f_{i-1}g_i=f_n
\end{equation}
Consider the left-hand side of (13)\\
$f_{i-1}g_i-d_{i-1}f_{i-2}g_i+f_{i-2}g_{i-1}\\
\hspace*{2cm}=f_{i-1}g_i-f_{i-2}(g_{i-1}+a_{i-1}b_{i-1}g_{i+1})+f_{i-2}g_{i-1}\\
\hspace*{2cm}=f_{i-1}g_i-a_{i-1}b_{i-1}f_{i-2}g_{i+1}\\
\hspace*{2cm}=c_{i-1}f_{i-2}e_ig_{i+1}-a_{i-1}b_{i-1}f_{i-2}g_{i+1}\\
\hspace*{2cm}=f_{i-2}g_{i+1}(c_{i-1}e_i-a_{i-1}b_{i-1})\\
\hspace*{2cm}=c_{i-1}f_{i-2}g_{i+1}(e_i-\frac{a_{i-1}b_{i-1}}{c_{i-1}})\\
\hspace*{2cm}=f_{i-1}g_{i+1}(e_i-d_i+c_i)$\\
Hence, we may write the formula (13) in the form:
\begin{equation}
    (c_i-d_i+e_i)f_{i-1}g_{i+1}=f_n
\end{equation}
Expanding the left-hand side of (15), yields\\
$$(c_if_{i-1})g_{i+1}-d_if_{i-1}g_{i+1}+f_{i-1}(e_ig_{i+1})=f_n$$
Therefore, by using (7) and (8), we get
$$f_ig_{i+1}-d_if_{i-1}g_{i+1}+f_{i-1}g_i=f_n$$
Consequently, the formula (15) is valid. This completes the proof.\\
Armed with the above result, we may now formulate the following result.
\\
\textbf{Theorem 1\cite{Huang1997}:} Let $A_n$ be a nonsingular tridiagonal matrix of the form (1). Then $B_n=A_n^{-1}=[\alpha_{ij}]_{1\le i,j\le n}$ is given by
\begin{equation}
  \alpha_{ij}=\left\{
  \begin{array}{lr} 
      \frac{1}{c_i-d_i+e_i}& \textit{if}\quad i=j \\
       -y_i\alpha_{i+1,j}& \textit{if}\quad i<j\\
       -v_{i-1}\alpha_{i-1,j}& \textit{if}\quad i>j
      \end{array}
      \right.\\
\end{equation}
where $y_i=\frac{a_i}{c_i},$ and $v_i=\frac{b_i}{e_{i+1}}, i=1, 2, \cdots, n-1$.\\
Repeated applications of (16) yields
\begin{equation}
  \alpha_{ij}=\left\{
  \begin{array}{lr} 
      \frac{1}{c_i-d_i+e_i}& \textit{if}\quad i=j \\
       (\prod_{r=i}^{j-1}\frac{p_r}{c_r})\alpha_{jj}& \textit{if}\quad i<j\\
       (\prod_{r=j}^{i-1}\frac{q_r}{e_{r+1}})\alpha_{jj}& \textit{if}\quad i>j
      \end{array}
      \right.\\
\end{equation}
\textbf{Theorem 2.} Consider the matrix $A_n$ in (1). If $det(A_n)\ne 0$ then the inverse of the matrix $A_n$ is given by $A_n^{-1}=[\alpha_{ij}]_{1\le i,j \le n}$ with $\alpha_{ij}$ given by:\\
\begin{equation}
\alpha_{ij}=\left\{
  \begin{array}{lr} 
      f_{i-1}g_{i+1}/f_n & \textit{if}\quad i=j \\
       (f_{i-1}g_{j+1}/f_n)(\prod_{r=i}^{j-1}p_r)& \textit{if}\quad i<j\\
       (f_{j-1}g_{i+1}/f_n)(\prod_{r=j}^{i-1}q_r)& \textit{if}\quad i>j
      \end{array}
      \right.
      \end{equation}
\textbf{Proof:} By using (7) and (8), we see that\\
$
        \alpha_{ii}=\frac{1}{c_i-d_i+e_i}=\frac{1}{\frac{f_i}{f_{i-1}}-d_i+\frac{g_i}{g_{i+1}}}=\frac{f_{i-1}g_{i+1}}{f_ig_{i+1}-d_if_{i-1}g_{i+1}+f_{i-1}g_i}=\frac{f_{i-1}g_{i+1}}{f_n}
$\\
by using \textbf{Lemma 1}.\\
Now, if $i<j$, then\\
$\alpha_{ij}=(\prod_{r=i}^{j-1}\frac{p_r}{c_r})\alpha_{jj}=\frac{\prod_{r=i}^{j-1}p_r}{\prod_{r=i}^{j-1}c_r}\frac{f_{j-1}g_{j+1}}{f_n}=\frac{\prod_{r=i}^{j-1}p_r}{\prod_{r=1}^{i-1}c_r\prod_{r=i}^{j-1}c_r}\frac{f_{j-1}g_{j+1}f_{i-1}}{f_n}\\
\hspace*{0.5cm}=\frac{f_{i-1}g_{j+1}}{f_n}\prod_{r=i}^{j-1}p_r$ having used (7).\\
Similarly, we can prove that:\\
For $i>j$, we have\\
$$\alpha_{ij}=\frac{f_{j-1}g_{i+1}}{f_n}\prod_{r=j}^{i-1}q_r$$\\
This completes the proof.\\
Following \cite{Kumar1993}, we see that\\
\begin{equation}
    \alpha_{ij}=[g_{j+1}(\sum_{k=1}^{j-1}\delta_{ik}f_{k-1}\prod_{r=k}^{j-1}p_r)+\delta_{ij}f_{j-1}g_{j+1}+f_{j-1}(\sum_{j+1}^{n}\delta_{ik}g_{k+1}\prod_{r=j}^{k-1}q_r)]/f_n
\end{equation}
is satisfied (see also \cite{Usmani1994}).\\
where $\delta_{ij}$ is the Kronecker delta symbol which is equal to $1$ if $i=j$ and $0$ if $i\ne j$. It should be noticed that we may obtain (18) from (19).\\
At this point, we can introduce the new numerical algorithm for inverting any general nonsingular tridiagonal matrix $A_n$ of the form (1).
\begin{algorithm}[!h]
 \caption{(To find the $n\times n$ inverse matrix of a general tridiagonal matrix\\ \hspace*{2.2cm} $A_n$ of the form (1)).}
  \textbf{Input:} The components $a_i, b_i, i=1, 2, \cdots, n-1$ and $d_i, i=1, 2, \cdots, n$.\\
  \textbf{Output:} The inverse matrix$A_n^{-1}=[\alpha_{ij}]_{1\le i,j\le n}$.\\
  \textbf{Step 1:} Compute $f_1, f_2, \cdots, f_n$ using the \textbf{DETGTRIF} algorithm. If $f_n=0$ then print"Singular matrix" and stop.\\
  \textbf{Step 2:} Compute $g_{n+1}, g_n, \cdots, g_1$ using the \textbf{DETGTRIG} algorithm.\\
  \textbf{Step 3:} Calculate the elements of the inverse matrix $\alpha_{ii}, i=1, 2, 3, \cdots, n$,\\
  $\alpha_{ij}, i=1, 2, 3, \cdots, n-1, j=i+1, i+2, \cdots, n$,\\
  and $\alpha_{ij}, j=1, 2, 3, \cdots, n-1, i=j+1, j+2, \cdots, n$ by using (18).\\
    \label{alg3}
  \end{algorithm}
  \\
The \textbf{Algorithm} \ref{alg3} will be referred to as \textbf{HINVGTRI} algorithm.\\
Using the \textbf{HINVGTRI} algorithm, a procedure has been written in Maple and can be obtained from the authors upon request.\\
\textbf{Remark 1.} If $a_k=0, 1\le k\le n-1$, then $\alpha_{ij}=0, i=1, 2, \cdots, k; j=k+1, k+2, \cdots, n$. Similarly, If $b_k=0, 1\le k\le n-1$, then $\alpha_{ij}=0, j=1, 2, \cdots, k; i=k+1, k+2, \cdots, n$.\\
\textbf{Remark 2.} If $f_k=0, 1\le k\le n-1$, then $\alpha_{k+1,j}=\alpha_{j,k+1}=0, j=k+1, k+2, \cdots, n$. Similarly, If $g_k=0, 2\le k\le n$, then $\alpha_{i,k-1}=\alpha_{k-1,i}=0, i=1, 2, \cdots, k-1$.\\
\textbf{Remark 3.} The inverse matrix $B_n$ of the matrix $A_n$ of the form in (1) can be written as the Hadamard product, also called element-wise product.
\begin{equation}
  B_n=R \circ S  
\end{equation}
where $R=(r_{ij})_{1\le i,j\le n}$ and $S=(s_{ij})_{1\le i,j\le n}$ are given by:
\begin{equation}
   r_{ij}=r_{ji}=\left\{
  \begin{array}{lr} 
      f_{i-1}g_{i+1}/\Delta & \textit{if}\quad i=j \\
       f_{i-1}g_{j+1}/\Delta& \textit{if}\quad i\ne j
      \end{array}
      \right. 
\end{equation}
\begin{equation}
   s_{ij}=\left\{
  \begin{array}{lr} 
      1 & \textit{if}\quad i=j \\
      \prod_{r=i}^{j-1}p_r& \textit{if}\quad i<j\\
      \prod_{r=j}^{i-1}q_r& \textit{if}\quad i>j
      \end{array}
      \right. 
\end{equation}
and 
\begin{equation}
    \alpha_{ij}=r_{ij}\times s_{ij}, i, j=1, 2, \cdots, n.
\end{equation}
Note that in formula (21), $\Delta= det(A_n)$. \\
\textbf{Remark 4.} For centrosymmetric tridiagonal matrices, we have $f_i=g_{n+1-i}, i=0, 1, 2, \cdots, n$. 

\section{Illustrative Examples}
In this section, We are going to present some illustrative examples.\\
\textbf{Example 1.} Consider $A_4=\left[
\begin{array}{cccc}
    2 & -1 & 0 & 0\\
   -2 & 2 & 1 & 0\\
   0 & 1 & 2 & 3\\
   0 & 0 & -1 & -3\\
 \end{array}
\right]$\\
Find $A_4^{-1}$.\\
\textbf{Solution}: We have\\
$a=[-1, 1, 3]^t$,
$b=[-2, 1, -1)]^t$, and
$d=[2, 2, 2, -3]^t$ by using the \textbf{DETGTRIF} algorithm, we get
$det(A_4)=0$. Therefore, $A_4^{-1}$ does not exist.\\
\textbf{Example 2.} Consider the $n\times n$ symmetric tridiagonal matrix $A_n$ given by \cite{Fonseca2007}:\\
$A_n=\left[
\begin{array}{ccccc}
    -1 & 1 &  & &\\
   1 & -2 & 1 & &\\
   \ddots & \ddots & \ddots & &\\
    &  & 1 & -2 & 1\\
    &  &  & 1 & -\frac{(n-1)}{n}\\
 \end{array}
\right]$\\
Find $A_n^{-1}$.\\
\textbf{Solution}: We have\\
$a=b=[1,1, \cdots, 1]^t$ and $d=[-1, -2, -2, \cdots, -2, -\frac{(n-1)}{n}]^t$,\\
\textbf{Step 1:} $c_i=-1, i=1, 2, \cdots, n-1$ and $c_n=\frac{1}{n}$. $e_i=-(\frac{i-1}{i}), i=n, n-1, \cdots, 2$\\ \hspace*{1.4cm} and $e_1=1$. $f_i=(-1)^i, i=1, 2, \cdots, n-1$ and $f_n=det(A_n)=\frac{(-1)^{n-1}}{n}$.\\ \hspace*{1.4cm}  $g_i=(-1)^{n+1-i}(\frac{i-1}{n}), i=n, n-1, \cdots, 2$ and $g_1=det(A_n)=\frac{(-1)^{n-1}}{n}$.\\
\textbf{Step 2:} When $i=j$ then $\alpha_{ii}=f_{i-1}g_{i+1}/det(A_n)=i$\\
\textbf{Step 3:} When $i<j$ then $\alpha_{ij}=(-1)^{j-i}f_{i-1}g_{j+1}/det(A_n)=j$\\
\textbf{Step 4:} When $i>j$ and since $A_n^{-1}$ is symmetric then $\alpha_{ij}=i$.\\
So $A_n^{-1}=[max(i,j)]_{1\le i,j\le n}$\\
\textbf{Example 3.} Compute $A^{-1}=\left[
\begin{array}{cccc}
    25 & -9 & 0 & 0\\
   -9 & 13 & -4 & 0\\
   0 & -4 & 5 & -1\\
   0 & 0 & -1 & 1\\
 \end{array}
\right]^{-1}$\\
\textbf{Solution}: We have\\
$a=b=[-9, -4, -1]^t$ and $d=[25, 13, 5, -1]^t$.\\
$f=[1, 25, 244, 820, 576]^t$,and
$g=[576, 36, 4, 1, 1]$.\\Consequently, $\Delta = 576$,\\
and\\
$A^{-1}=\frac{1}{576}\left[
\begin{array}{cccc}
    36 & 36 & 36 & 36\\
   36 & 100 & 100 & 100\\
   36 & 100 & 244 & 244\\
   36 & 100 & 244 & 820\\
 \end{array}
\right]$\\
\textbf{Example 4.} Consider the centrosymmetric matrix $A_n=\left[
\begin{array}{ccccc}
    2 & -1 &  & &\\
   -1 & 2 & -1 & &\\
   \ddots & \ddots & \ddots & &\\
    &  & -1 & 2 & -1\\
    &  &  & -1 & 2\\
 \end{array}
\right]$\\
Find $A^{-1}$.\\
\textbf{Solution}: We have $a=b=[-1, -1, \cdots, -1]^t$ and $d=[2, 2, \cdots, 2]^t$.\\
Since $A_n$ is centrosymmetric matrix, then $f_i=i+1, g_i=n+2-i, i=1, 2, \cdots, n$.\\
$A^{-1}=X \circ R, \quad x_{ij}=1 for 1\le i,j\le n$,\\
Thus $A^{-1}=R, r_{ij}=r_{ji}=\frac{1}{n+1}f_{i-1} g_{j+1}=\frac{1}{n+1}(i)(n+1-j)=\frac{i(n+1-j)}{n+1}$.\\\\
\textbf{Example 5.} Consider $A_5=\left[
\begin{array}{ccccc}
    1 & -1 & 0 & 0 & 0\\
   -1 & 3 & -1 & 0 & 0\\
   0 & -1 & 1 & -1 & 0\\
   0 & 0 & -1 & 1 & -1\\
   0 & 0 & 0 & -1 & 1\\
 \end{array}
\right]$\\
Find $A_5^{-1}$.\\
\textbf{Solution}: We have $a=b=[-1, -1, -1, -1]^t$ and 
$d=[1, 3, 1, 1, 1]^t$.\\
Thus, $f=[1, 1, 2, 1, -1, -2]^t$, $g=[-2, -3, -1, 0, 1, 1]^t$ and\\\\
$A^{-1}=\left[
\begin{array}{ccccc}
    \frac{3}{2} & \frac{1}{2} & 0 & -\frac{1}{2} & -\frac{1}{2}\\
  \frac{1}{2} & \frac{1}{2} & 0 & -\frac{1}{2} & -\frac{1}{2}\\
  0 & 0 & 0 & -1 & -1\\
  -\frac{1}{2} & -\frac{1}{2} & -1 & -1 & -\frac{1}{2}\\
  -\frac{1}{2} & -\frac{1}{2} & 1 & -1 & \frac{1}{2}\\
 \end{array}
\right]$.\\
\bibliography{mybibfile}

\end{document}